\documentstyle{article}

\newcommand{\eL}{\mbox{$\cal L$}}
\newcommand{\Pe}{\mbox{$\cal P$}}

\newcommand{\kon}{\wedge}

\newcommand{\str}{\rightarrow}
\newcommand{\rts}{\leftarrow}
\newcommand{\mj}{\mbox{\bf 1}}

\newcommand{\df}{\mbox{\scriptsize{\it df}}}
\newcommand{\HDS}{\vrule width0pt height2.3ex depth1.05ex\displaystyle}

\newcommand{\Bx}{{\raisebox{-1pt}{
\begin{picture}(8,8)
\put(3,2){\makebox(0,0){$\Box$}}
\put(3,3.5){\makebox(0,0){$\times$}}
\end{picture}}}\,}

\newcommand{\Bk}{{\raisebox{-1pt}{
\begin{picture}(8,8)
\put(3,2){\makebox(0,0){$\Box$}}
\put(2.9,3.3){\makebox(0,0){\scriptsize +}}
\end{picture}}}\,}

\def\cirk{\,{\raisebox{.3ex}{\tiny $\circ$}}\,}
\def\ks{\mbox{\footnotesize$\;\xi\;$}}
\def\b#1#2{\stackrel{\raisebox{-2pt}{\mbox{\tiny $#1$}}}
{\raisebox{0pt}{$b$}}^{\raisebox{-7pt}{\scriptsize $#2$}}}

\def\c#1{\stackrel{\raisebox{-2pt}{\mbox{\tiny $\,#1$}}}
{\raisebox{0pt}{$c$}}}

\def\d#1#2{\stackrel{\raisebox{-2pt}{\mbox{\tiny $\,#1$}}}
{\raisebox{0pt}{$\delta$}}^{\raisebox{-7pt}{\scriptsize $#2$}}}

\def\w#1{\stackrel{\raisebox{-2pt}{\mbox{\tiny $\,#1$}}}
{\raisebox{0pt}{$w$}}}

\def\f#1#2{{{\HDS #1}\over{\HDS #2}}}

\def\set{\mbox{\it Set}}

\def\pl{\!+\!}
\def\mn{\!-\!}

\def\ckm{\c{\kon}^{\raisebox{-4pt}{\scriptsize $m$}}}

\newcommand{\SMC}{\mbox{\bf SMC}}

\begin{document}

\title{Relevant Categories and Partial Functions}
\author{{\sc Kosta Do\v sen} and {\sc Zoran Petri\' c}\\[0.5cm]
Mathematical Institute, SANU \\
Knez Mihailova 35, p.f. 367 \\
11001 Belgrade, Serbia \\
email: \{kosta, zpetric\}@mi.sanu.ac.yu}
\date{}
\maketitle

\begin{abstract}
\noindent A relevant category is a symmetric monoidal closed
category with a diagonal natural transformation that satisfies
some coherence conditions. Every cartesian closed category is a
relevant category in this sense. The denomination \emph{relevant}
comes from the connection with relevant logic. It is shown that
the category of sets with partial functions, which is isomorphic
to the category of pointed sets, is a category that is relevant,
but not cartesian closed.
\end{abstract}

\vspace{.3cm}


\vspace{.5ex}

\noindent {\it Keywords$\,$}: symmetric monoidal closed
categories, diagonal natural transformation, intuitionistic
relevant logic, partial functions, pointed sets

\vspace{0.5cm}

\begin{flushright}
{\it To the memory of Aleksandar Kron}
\end{flushright}

\vspace{0.5cm}

\baselineskip=1.2\baselineskip

\section{Introduction}

It is well known that the category of pointed sets with
point-preserving functions as arrows is a symmetric monoidal
closed category (see \cite{EK66a}, Section IV.1). The symmetric
monoidal closed structure in this category is provided by the
smash product and internal hom-sets (see \S 3 below). This
category, which is isomorphic to the category of sets with partial
functions, has however a richer structure than just symmetric
monoidal closed.

If to the assumptions for symmetric monoidal categories we add a
diagonal natural transformation with appropriate equations between
arrows, then one obtains a notion of monoidal category for which a
coherence theorem is proved in \cite{P02} with respect to
relations on finite ordinals. We will call these categories
\emph{relevant monoidal categories}, because the types of arrows
in the relevant monoidal category freely generated by a set of
propositional letters corresponds to sequents in the
multiplicative (i.e.\ intensional) conjunction (i.e.\ fusion)
fragment of relevant logic, including the multiplicative constant
true proposition $\top$. This fragment, which is the same both in
intuitionistic and in classical versions of relevant logic,
catches essentially the structural rules of relevant logic, on
which the whole structure of this logic rests. The equations for
relevant monoidal categories stem from \cite{DP96}. (An incomplete
set of these equations may be found in \cite{Jac94}, Definition
2.1(i).)

Symmetric monoidal closed categories that are also relevant
monoidal with respect to the same monoidal structure will be
called \emph{relevant monoidal closed} categories. The relevant
monoidal closed category \textbf{RMC} freely generated by a set of
propositional letters corresponds to the multiplicative
conjunction-$\top$-implication fragment (both intuitionistic and
classical) of relevant logic. We have that $A$ is the source and
$B$ the target of an arrow of \textbf{RMC} iff $A\str B$ is a
theorem of the multiplicative conjunction-$\top$-implication
fragment of the relevant logic \textbf{R} (see \cite{AB75}, where
multiplicative conjunction, i.e.\ fusion, is called
\emph{co-tenability}).

The category of pointed sets has a diagonal natural transformation
with respect to the smash product, which makes of it a relevant
monoidal closed category. Moreover, it has also finite products
and coproducts (including the empty ones), where product is
different from the smash product. With this additional structure,
we obtain also operations corresponding to the additive (lattice)
connectives of relevant logic, without distribution of additive
conjunction over additive disjunction. The relevant logic
\textbf{R} has this distribution, but a version of \textbf{R}
without it also exists. (Linear logic lacks this distribution.) We
are still within a fragment of relevant logic common to its
intuitionistic and classical versions. This fragment catches
presumably the whole positive fragment of intuitionistic relevant
logic.

The connection between relevant logic and the category of pointed
sets is reminiscent of the connection that exists between
intuitionistic logic and the category \emph{Set} of sets with
functions. As intuitionistic propositional logic may be identified
with the bicartesian closed category (see \cite{LS86}, Section
I.8) freely generated by a set of propositional letters, so the
positive fragment of intuitionistic relevant logic may be
identified with a free relevant category such as we will
introduce. And as \emph{Set} is the prime example of a bicartesian
closed category, so the category of pointed sets may be the prime
example of a relevant category.

Inspired by some ideas of Belnap, which are derived from Scott's
models for the untyped lambda calculus, Helman found in \cite{H77}
that typed lambda terms in beta-normal form that code proofs in
the additive conjunction-implication fragment of the relevant
logic \textbf{R} can be interpreted in the hierarchy of pointed
sets with product and sets of point-preserving functions. A
connection between relevant logic and the category of pointed sets
was also investigated by Szabo in \cite{Sz83}, but with an
approach different from ours---in particular as far as
distribution of product over coproduct is concerned. It was
prefigured by Jacobs in \cite{Jac94} (Example 2.3(i)) that the
category of pointed sets is a relevant monoidal category, though
the notion of relevant monoidal category of that paper differs
from ours.

\section{Relevant categories}

The objects of the category \textbf{SyMon} are the formulae of the
propositional language ${\eL_{\top,\kon}}$, generated from a set
\Pe\ of propositional letters with the nullary connective, i.e.\
propositional constant, $\top$ and the binary connective $\kon$.
We use $p,q,r,\ldots\,$, sometimes with indices, for propositional
letters, and $A,B,C,\ldots\,$, sometimes with indices, for
formulae. As usual, we omit the outermost parentheses of formulae
and other expressions later on.

To define the arrows of \textbf{SyMon}, we define first
inductively a set of expressions called the {\it arrow terms}.
Every arrow term of \textbf{SyMon} will have a {\it type}, which
is an ordered pair of formulae of ${\eL_{\top,\kon}}$. We write
${f\!:A\vdash B}$ when the arrow term $f$ is of type ${(A,B)}$.
(We use the turnstile $\vdash$ instead of the more usual $\str$,
which we reserve for a connective and a biendofunctor.) We use
$f,g,h,\ldots\,$, sometimes with indices, for arrow terms.

For all formulae $A$, $B$ and $C$ of ${\eL_{\top,\kon}}$ the
following {\it primitive arrow terms}:

\begin{tabbing}
\centerline{$\mj_A\!: A\vdash A$,}
\\[1.5ex]
\centerline{$\b{\kon}{\str}_{A,B,C}\::A\kon(B\kon C)\vdash (A\kon
B)\kon C$,\quad $\b{\kon}{\rts}_{A,B,C}\::(A\kon B)\kon C\vdash
A\kon (B\kon C)$,}
\\[1.5ex]
\centerline{$\c{\kon}_{A,B}\::A\kon B\vdash B\kon A$,}
\\[1.5ex]
\centerline{$\d{\kon}{\str}_A: A\kon\top \vdash A$,\quad
$\d{\kon}{\rts}_A: A \vdash A\kon\top$}

\end{tabbing}

\noindent are arrow terms of \textbf{SyMon}. If ${g\!:A\vdash B}$
and ${f\!:B\vdash C}$ are arrow terms of \textbf{SyMon}, then
${f\cirk g\!:A\vdash C}$\index{composition@composition $\cirk$} is
an arrow term of \textbf{SyMon}; and if ${f\!:A\vdash D}$ and
${g\!:B\vdash E}$ are arrow terms of \textbf{SyMon}, then ${f\kon
g\!:A\kon B\vdash D\kon E}$ is an arrow term of \textbf{SyMon}.
This concludes the definition of the arrow terms of
\textbf{SyMon}.

Next we define inductively the set of {\it equations} of
\textbf{SyMon}, which are expressions of the form ${f=g}$, where
$f$ and $g$ are arrow terms of \textbf{SyMon} of the same type. We
stipulate first that all instances of ${f=f}$ and of the following
equations are equations of \textbf{SyMon}:

\begin{tabbing}
\mbox{\hspace{1em}}\= $({\mbox{{\it cat}~1}})$\quad\quad\= $f\cirk
\mj_A=\mj_B\cirk f=f\!:A\vdash B$,\index{cat1@({\mbox{{\it
cat}~1}}) equation}
\\*[1ex]
\> $({\mbox{{\it cat}~2}})$\> $h\cirk (g\cirk f)=(h\cirk g)\cirk
f$,\index{cat2@({\mbox{{\it cat}~2}}) equation}
\\[1.5ex]
 \> $(\kon\, 1)$\> $\mj_A\kon\mj_B=\mj_{A\kon B}$,
\\*[1ex]
\> $(\kon\, 2)$\> $(g_1\cirk f_1)\kon(g_2\cirk f_2)=(g_1\kon
g_2)\cirk(f_1\kon f_2)$,
\\[1.5ex]
for $f\!:A\vdash D$, $g\!:B\vdash E$ and $h\!:C\vdash F$,
\\*[1ex]
\> $\mbox{($\b{\kon}{\str}$ {\it nat})}$\>  $((f\kon g)\kon
h)\cirk\!\b{\kon}{\str}_{A,B,C}\:=\:
\b{\kon}{\str}_{D,E,F}\!\cirk(f\kon (g\kon h))$,
\\[1.5ex]
\> $\mbox{($\c{\kon}$ {\it nat})}$\> $(g\kon
f)\cirk\!\c{\kon}_{A,B}\:=\:\c{\kon}_{D,E}\!\cirk(f\kon g)$,
\\[1.5ex]
\> ($\d{\kon}{\str}$~{\it nat})\> $f\cirk\!
\d{\kon}{\str}_A\:=\:\d{\kon}{\str}_B\!\! \cirk(f\kon\mj_\top)$,
\\[1.5ex]
\>$(\b{\kon}{}\b{\kon}{})$\quad\=
$\b{\kon}{\str}_{A,B,C}\!\cirk\!\b{\kon}{\rts}_{A,B,C}\;
=\mj_{(A\kon B)\kon C}$,\quad\quad\=
$\b{\kon}{\rts}_{A,B,C}\!\cirk\!\b{\kon}{\str}_{A,B,C}\;=\mj_{A\kon(B\kon
C)}$,
\\[1.5ex]
\>$(\b{\kon}{}\!5)$\> $\b{\kon}{\str}_{A\kon B, C,D}\!\cirk\!
\b{\kon}{\str}_{A,B,C\kon D}\;=(\b{\kon}{\str}_{A,B,C}\kon
\:\mj_D)\cirk \!\b{\kon}{\str}_{A,B\kon
C,D}\!\cirk(\mj_A\:\kon\b{\kon}{\str}_{B,C,D})$,
\\[1.5ex]
\>$(\c{\kon}\c{\kon})$\>
$\c{\kon}_{B,A}\!\cirk\!\c{\kon}_{A,B}\;=\mj_{A\kon B}$,
\\[1ex]
\>$(\b{\kon}{}\c{\kon})$\> $\c{\kon}_{A,B\kon
C}\;=\;\b{\kon}{\str}_{B,C,A}\!\cirk(\mj_B\:\kon\c{\kon}_{A,C})
\cirk\!\b{\kon}{\rts}_{B,A,C}\!\cirk(\c{\kon}_{A,B}\kon\:
\mj_C)\cirk \!\b{\kon}{\str}_{A,B,C}$,
\\[1.5ex]
\>($\d{\kon}{}\d{\kon}{}$)\>
$\d{\kon}{\str}_A\!\cirk\!\d{\kon}{\rts}_A\:=\mj_A$, \quad\quad
$\d{\kon}{\rts}_A\!\cirk\!\d{\kon}{\str}_A\:= \mj_{A\kon\top}$,
\\*[1ex]
\>($\b{\kon}{}\d{\kon}{}$)\>
$\b{\kon}{\str}_{A,B,\top}\;=\;\d{\kon}{\rts}_{A\kon B}\!\cirk
(\mj_A\:\kon \d{\kon}{\str}_B)$.

\end{tabbing}

The set of equations of \textbf{SyMon} is closed under symmetry
and transitivity of equality and under the rules

\[
(\mbox{\it cong~}\ks\!)\quad \f{f=f_1 \quad \quad \quad g=g_1}
{f\ks g=f_1\ks g_1}\index{congxi@(\mbox{\it cong~}\ks) rule}
\]

\noindent where ${\!\ks\!\in\{\cirk,\kon\}}$; if $\!\ks\!$ is
$\cirk$, then ${f\cirk g}$ is defined (namely, $f$ and $g$ have
appropriate, composable, types), and analogously for ${f_1\cirk
g_1}$. This concludes the definition of the equations of
\textbf{SyMon}.

On the arrow terms of \textbf{SyMon} we impose the equations of
\textbf{SyMon}. This means that an arrow of \textbf{SyMon} is an
equivalence class of arrow terms of \textbf{SyMon} defined with
respect to the smallest equivalence relation such that the
equations of \textbf{SyMon} are satisfied (see \cite{DP04},
Section 2.3, for details).

The equations $(\kon\, 1)$ and $(\kon\, 2)$ say that $\kon$ is a
biendofunctor (a 2-endofunctor, in the terminology of \cite{DP04},
Section 2.4). Equations with ``\emph{nat}'' in their names, like
those in the list above, say that $\b{\kon}{\str}$, $\c{\kon}$,
etc.\ are natural transformations.

The category \textbf{SyMon} is the free symmetric monoidal
category in the sense of \cite{ML71} (Chapter VII) generated by
the set \Pe.

The category \textbf{ReMon} is defined as the category
\textbf{SyMon} with the following additions. We have the
additional primitive arrow terms

\[
\w{\kon}_A:A\vdash A\kon A,
\]

\noindent and the following additional equations:

\begin{tabbing}
\mbox{\hspace{5em}}\= (${\d{\kon}{\str}}$~{\it
nat})\mbox{\hspace{1.8em}} \= \kill

\> \mbox{($\w{\kon}$ {\it nat})} \> $(f\kon
f)\cirk\!\w{\kon}_A\:=\:\w{\kon}_D\!\cirk f$,
\\[1.5ex]
\> ($\b{\kon}{}\w{\kon}$) \>
$\b{\kon}{\str}_{A,A,A}\!\cirk(\mj_A\:\kon
\w{\kon}_A)\cirk\!\w{\kon}_A\: =\,(\w{\kon}_A\kon
\:\mj_A)\cirk\!\w{\kon}_A$,
\\*[1.5ex]
\> ($\c{\kon}\w{\kon}$) \>
$\c{\kon}_{A,A}\!\cirk\!\w{\kon}_A\;=\;\w{\kon}_A$,
\\[1.5ex]
for $\ckm_{A,B,C,D}\;=_{\df}\; \b{\kon}{\str}_{A,C,B\kon
D}\!\!\cirk (\mj_A\kon(\b{\kon}{\rts}_{C,B,D}\!\!\cirk
(\c{\kon}_{B,C}\!\kon\mj_D)\cirk\!\!\b{\kon}{\str}_{B,C,D}))
\cirk\!\!\b{\kon}{\rts}_{A,B,C\kon D}\::$
\\*[.5ex]
\` $(A\kon B)\kon(C\kon D)\vdash (A\kon C)\kon (B\kon D)$,
\\*[1ex]
\> ($\b{\kon}{}\c{\kon}\w{\kon}$) \> $\w{\kon}_{A\kon
B}\;=\;\ckm_{A,A,B,B}\!\cirk (\w{\kon}_A\kon \w{\kon}_B)$,
\\[1.5ex]
\> \mbox{($\w{\kon}\d{\kon}{}$)} \> $\w{\kon}_{\top}\:=\;
\d{\kon}{\rts}_{\top}$.

\end{tabbing}

\noindent (These equations may be found in \cite{DP96}; they are
also in \cite{Jac94}, Definition 2.1(i), but with
($\b{\kon}{}\c{\kon}\w{\kon}$) lacking.)

A \emph{relevant monoidal category} is a symmetric monoidal
category that has in addition a natural transformation $\w{\kon}$
that satisfies the equations of \textbf{ReMon}. The category
\textbf{ReMon} is the free relevant monoidal category generated by
the set \Pe. A coherence theorem is proved for this category in
\cite{P02} with respect to the category whose arrows are relations
between finite ordinals. This means that there is a faithful
functor from \textbf{ReMon} into the latter category.

The category \SMC\ is defined as the category \textbf{SyMon} with
the following additions. We have an additional binary connective
$\str$, and the additional primitive arrow terms

\[
\varepsilon_{A,B}\!: A\kon (A\str B)\vdash B, \quad\quad
\eta_{A,B}\!: B\vdash A\str (A\kon B);
\]

\noindent on arrow terms we have the additional unary operations
${A\str}$, for every object $A$, such that for ${f\!:B\vdash C}$
we have the arrow term ${A\str f\!:A\str B\vdash A\str C}$.

The equations of \SMC\ are obtained by assuming the following
additional equations:

\begin{tabbing}
\mbox{\hspace{5em}} \= (${\d{\kon}{\str}}$~{\it
nat})\mbox{\hspace{1.8em}} \= $f\cirk
\d{\kon}{\str}_A=\d{\kon}{\str}_B\cirk(f\kon\mj_\top)$,\kill

\> ($A\str$~1)\> $A\str \mj_B = \mj_{A\str
B}$,\index{Aarrow1@($A\str$~1) equation}
\\*[1ex]
\> ($A\str$~2)\> $A\str (f\cirk g) = (A\str f)\cirk (A\str
g)$,\index{Aarrow2@($A\str$~2) equation}
\\[1.5ex]
\> ($\varepsilon$~{\it nat})\>
$f\cirk\varepsilon_{A,B}=\varepsilon_{A,C}\cirk(\mj_A\kon(A\str
f))$,
\\*[1ex]
\> ($\eta$~{\it nat})\> $(A\str(\mj_A\kon
f))\cirk\eta_{A,B}=\eta_{A,C}\cirk f$,
\\[1.5ex]
\>($\varepsilon\eta\;\kon$)\> $\varepsilon_{A,A\kon
B}\cirk(\mj_A\kon\eta_{A,B})\:=\mj_{A\kon B}$,
\\*[1ex]
\> ($\varepsilon\eta\str$)\>
$(A\str\varepsilon_{A,B})\cirk\eta_{A,A\str B}=\mj_{A\str B}$,
\end{tabbing}

\noindent and the following additional rule:

\[
\f{f=g}{A\str f=A\str g}
\]

The equations (${A\str}$~1) and (${A\str}$~2) say that $A\str$ is
a functor, while ($\varepsilon\eta\;\kon$) and
($\varepsilon\eta\str$) are the triangular equations of an
adjunction (see \cite{ML71}, Section IV.1). The category \SMC\ is
the free symmetric monoidal closed category generated by the set
\Pe\ (see \cite{ML71}, Section VII.7).

The category \textbf{RMC} is defined by combining the definitions
of \textbf{ReMon} and \SMC. \emph{Relevant monoidal closed
categories} are symmetric monoidal closed categories that are also
relevant monoidal with respect to the same monoidal structure. The
category \textbf{RMC} is the free relevant monoidal closed
category generated by the set \Pe.

A \emph{positive intuitionistic relevant} category is a relevant
monoidal closed category that has all finite products and
coproducts (including the empty ones).

\section{The category of pointed sets}

\noindent {\large\bf 3.1}\quad The category of pointed sets
$\set_{\ast}$ is the category whose objects are sets with a
distinguished element $*$, and whose arrows are functions $f$ such
that $f(*)=*$. This category is isomorphic to the category of sets
with partial functions, i.e.\ relations that are single-valued,
but not necessarily defined on the whole domain.

We have the following special objects and operations on objects in
$\set_{\ast}$:

\begin{tabbing}
\quad\quad\quad\quad I$\;=\{\ast\}$, \quad\quad \=
$a'=\{(x,\ast)\mid x\in a-{\mbox{\rm I}}\}$, \quad\quad
$b''=\{(\ast,y)\mid y\in b-{\mbox{\rm I}}\}$,
\\[1.5ex]
\> $a\otimes b\:$ \= = \= $((a-{\mbox{\rm I}})\times(b-{\mbox{\rm
I}}))\cup {\mbox{\rm I}}$,
\\[.5ex]
\> $a\Bx b\:$ \> = \> $(a\otimes b)\cup a'\cup b''$,
\\[.5ex]
\> $a\Bk b\:$ \> = \> $a'\cup b''\cup {\mbox{\rm I}}$.
\end{tabbing}

\noindent Let $\top$ in $\set_{\ast}$ be ${\mbox{\rm I}}\cup
\{x\}$, where $x\neq\ast$, and let $\kon$ be $\otimes$, which is
the \emph{smash} product. For $a$ and $b$ objects of
$\set_{\ast}$, let $a\str b$ be the union of I with the set of
arrows of $\set_{\ast}$ from $a$ to $b$ without the arrow with
constant value $*$. It is well known that with these operations on
objects $\set_{\ast}$ is a symmetric monoidal closed category (see
\cite{EK66a}, Section IV.1; see \cite{DP97}, Section 6, for more
details). The importance of $\set_{\ast}$ for checking equations
between arrows in \SMC\ is demonstrated by Soloviev in \cite{S97}.

It can also be easily shown that with ${\w{\kon}_a\!\!(x)}$ being
$(x,x)$ for $x\neq\ast$, and $\ast$ otherwise, $\set_{\ast}$ is a
relevant monoidal closed category. It was prefigured in
\cite{Jac94} (Example 2.3(i)) that $\set_{\ast}$ is a relevant
monoidal category, though the notion of relevant monoidal category
of that paper differs from ours, as explained in the preceding
section.

We also have arrows of the type of projections for the smash
product, which are defined in an obvious way, but these arrows do
not make natural transformations. The category $\set_{\ast}$ is
not a cartesian category with the smash product. But $\set_{\ast}$
is a bicartesian category (i.e.\ a category with all finite
products and coproducts) with binary product being $\Bx$ (which
corresponds to cartesian product) and binary coproduct being
$\Bk\!\!$, while I is both the terminal and the initial object,
i.e.\ the empty product and coproduct. It is shown in \cite{DP04}
(Sections 9.7, 12.4 and 13.4) that $\set_{\ast}$ is a bicartesian
category of a particular kind, called there
\emph{zero-dicartesian}.

So $\set_{\ast}$ is a positive intuitionistic relevant category in
the sense of the preceding section.

\noindent {\large \bf 3.2}\quad In \cite{DP97} one can find a
characterization of the objects isomorphic in \textbf{SMC}, and
hence also in every symmetric monoidal closed category. Two
objects $A$ and $B$ of \textbf{SMC} are isomorphic iff one can
derive $A=B$ in the equational calculus \textbf{S} whose axioms
are the axioms of commutative monoids with respect to $\kon$ and
$\top$ and the following two equations:

\begin{tabbing}
\hspace{12em}\=$\top\str C$ \= $=\:C$,
\\
\>\>$(A\kon B)\str C\!$\'$=\:B\str(A\str C)$.
\end{tabbing}

The proof of this equivalence in \cite{DP97} is based on a proof
of an analogous equivalence where \textbf{SMC} is replaced by
{\emph{Fin}$\set_*$}, which is the category of finite pointed sets
with point-preserving functions (a full subcategory of $\set_*$),
and $A$ and $B$ are \emph{diversified}, which means that no
propositional letter occurs more than once in these formulae. The
equational calculus \textbf{S} axiomatizes all the equations
between diversified formulae that hold in natural numbers, where
formulae are understood as arithmetical terms such that
propositional letters are variables ranging over natural numbers,
$\top$ is 1, the operation $\kon$ is multiplication, and $m\str n$
is $(n\pl 1)^m\mn 1$.

We conjecture that these two assertions concerning \textbf{S}, its
arithmetical interpretation and {\emph{Fin}$\set_*$} are also true
when the restriction to diversified formulae is lifted. If this
conjecture were true, then the calculus \textbf{S} would
characterize not only all the formulae isomorphic in \textbf{SMC},
but also in \textbf{RMC}.

\baselineskip=0.84\baselineskip

\end{document}